\documentclass[11pt]{article}
\usepackage{IEEEtrantools}
\usepackage{mathtools, nccmath}
\usepackage{latexsym}
\usepackage{amsfonts}
\usepackage{amssymb}
\usepackage{eqnarray, amsmath,amsfonts,amsthm}
\usepackage{mathrsfs}
\usepackage{enumerate}
\usepackage{dsfont}
\usepackage{hyperref}
\usepackage{graphicx}
\usepackage{epstopdf}
\usepackage{tikz}

\usetikzlibrary{shapes,arrows.meta,positioning}
\usepackage{cleveref}
\usepackage{tikz-cd}
\usepackage{psfrag}

\usepackage{multirow}
\usepackage{tikz}
\usetikzlibrary{arrows,decorations.markings}
\usepackage{bookmark}
\usepackage{hyperref}
\hypersetup{
     colorlinks   = true,
     citecolor    = blue
}

\newcommand{\newsection}[1]{\setcounter{equation}{0}
\setcounter{dfn}{0}
\section{#1}}

\newtheorem{dfn}{Definition}[section]
\newtheorem{thm}[dfn]{Theorem}
\newtheorem{lmma}[dfn]{Lemma}
\newtheorem{ppsn}[dfn]{Proposition}

\newtheorem{rmrk}[dfn]{Remark}
\newtheorem{notation}[dfn]{Notation}
\newtheorem*{theoremA}{Theorem A}
\newtheorem*{theoremB}{Theorem B}

\DeclareMathOperator*{\dprime}{\prime \prime}

\newcommand{\ot}{\otimes}

\newcommand{\dt}{\Delta_n}
\newcommand{\bbc}{\mathbb{C}}
\newcommand{\bbz}{\mathbb{Z}}
\newcommand{\bbn}{\mathbb{N}}

\def \qed { \mbox{}\hfill
$\Box$\vspace{1ex}}

\title{Sections and Chapters}
\begin{document}

\tikzset{->-/.style={decoration={
  markings,
  mark=at position #1 with {\arrow{>}}},postaction={decorate}}}
  \tikzset{-<-/.style={decoration={
  markings,
  mark=at position #1 with {\arrow{<}}},postaction={decorate}}}

\author{\sc{Keshab Chandra Bakshi,\,\,Satyajit Guin,\,\,Guruprasad}}
\title{Intersection of conjugate Hadamard subfactors arising from Fourier matrices }
\maketitle


\begin{abstract}

Given two distinct complex Hadamard matrices belonging to the same equivalence class generated by the tensor products of Fourier matrices, we show that if the corresponding Hadamard subfactors are conjugate, then their intersection is a factor with finite Jones index. We compute the index of the intersection explicitly and determine its relative commutant. Furthermore, we precisely characterize when these intersections give rise to vertex model subfactors, thereby extending our earlier results in low dimensions. As an application, we derive an explicit formula for the Connes-St{\o}rmer relative entropy associated with these intersections. These results reveal how the internal algebraic structure of complex Hadamard matrices governs the relative position and entropic behaviour of the subfactors.
\end{abstract}
\bigskip

{\bf AMS Subject Classification No.:} {\large 46}L{\large 37}; {\large 46}L{\large 55}
\smallskip


\hypersetup{linkcolor=blue}
\bigskip

\newsection{Introduction}\label{Sec1}

The study of {\em two subfactors} inside a common ambient factor reveals subtle structural and categorical phenomena that go far beyond the single-subfactor theory \cite{Jo}. Given two subfactors $N_1, N_2 \subset M$ of a type $\mathrm{II}_1$ factor, the intersection $N_1\cap N_2$ can display a rich and sometimes unexpected behavior. In general, $N_1 \cap N_2$ need not be a factor, and even when it is, its Jones index in $M$ may be infinite. For instance, Jones constructed (in \cite{Jo1}, see also \cite{BG2}) an example of $N_1, N_2 \subset M$ with $[M:N_1] = [M:N_2] = 2$ for which the (Pimsner--Popa) index of $N_1 \cap N_2$ in $M$ is infinite. Latter, Jones and Xu established in \cite{JX} that the intersection of two finite-index subfactors has finite index if and only if the corresponding angle operator \cite{SW} has finite spectrum.

In the absence of a general theory describing intersections of subfactors, it becomes natural to examine specific but meaningful families where explicit computation is possible. Motivated by this, we focus on {\em Hadamard} (or {\em spin model}) subfactors, a class introduced and emphasized by Jones, which lie at the crossroads of operator algebras, combinatorics, and quantum symmetries. These subfactors, arising from complex Hadamard matrices, offer a concrete setting where one can analyse the behaviour of intersections in depth (see \cite{BG1, BG2, BGG, BGG2} for recent developments).

The main goal of this paper is to investigate the intersection of two \emph{conjugate Hadamard subfactors} $R_U, R_V \subset R$ arising from distinct complex Hadamard matrices $U, V$ belonging to the same equivalence class generated by tensor products of Fourier matrices. We show that the intersection $R_U \cap R_V$ is always a factor of finite Jones index, and we determine this index explicitly in terms of an associated finite abelian group. Moreover, we compute the relative commutant of the intersection and identify the precise condition under which it yields a {\em vertex model subfactor}.

Our first main result asserts that factoriality of the intersection holds whenever the Hadamard subfactors are conjugate.

\begin{theoremA}
Let $U$ and $V$ be two distinct complex Hadamard matrices in the Hadamard equivalence class generated by tensor products of Fourier matrices $F_{n_1}\otimes\ldots\otimes F_{n_k}$ with $n_1,\ldots,n_k \in \bbn$. If the Hadamard subfactors $R_U\subset R$ and $R_V\subset R$ are conjugate to each other, then their intersection $R_U\cap R_V$ is a factor.
\end{theoremA}
We remark that in Theorem $5.2$, \cite{BGG2}, we  have provided a complete characterization exactly when the Hadamard subfactors are conjugate.

We further compute the Jones index and the relative commutant of the intersection.

\begin{theoremB}
Let $U,V\in [F_{n_1} \otimes \cdots \otimes F_{n_k}]$ be distinct complex Hadamard matrices with associated Hadamard subfactors $R_U, R_V \subset R$. If $R_U\subset R$ and $R_V\subset R$ are conjugate, then
\[
[R : R_U \cap R_V] = \frac{N^2}{|H|}\,,
\]
where $N=n_1\ldots n_k$ and $H$ is a subgroup of $\mathbb{Z}_{n_1} \times \cdots \times \mathbb{Z}_{n_k}$ with
\[
|H| = \dim\big(\mathrm{Ad}_U(\Delta_N) \cap \mathrm{Ad}_V(\Delta_N)\big).
\]
Here $\Delta_N$ is the abelian subalgebra of diagonal matrices in $M_N(\mathbb{C})$. Moreover,
\[
(R_U \cap R_V)' \cap R = (\mathrm{Ad}_U(\Delta_N) \cap \mathrm{Ad}_V(\Delta_N))' \cap \mathbb{C}^N.
\]
\end{theoremB}

We also prove that every possible value of the index of the form $N^2/|H|$ can indeed be realized by a suitable pair $(U, V)$ of Hadamard matrices in the same equivalence class. This shows that although each Hadamard subfactor $R_U \subset R$ has the same Jones index $N$, the index of their intersection varies with the specific algebraic form of $U$ and $V$, revealing a delicate dependence on the underlying matrix data. We also obtain a precise characterization of exactly when the subfactors $R_U\cap R_V\subset R$ correspond to {\em vertex model} subfactors, and it generalizes aspects of our earlier work in low dimensions \cite{BG1,BGG}. Since vertex models are connected to fusion categories, quantum groups, and statistical mechanics, identifying when a subfactor arises from such a model provides valuable information for classification and theoretical understanding.

Finally, we analyse the {\em Connes--St{\o}rmer relative entropy} associated with these intersections and obtain explicit bounds and formulae that quantify their mutual position within $R$. In particular, we derive an exact expression for the modified relative entropy in terms of the entries of $U^*V$, demonstrating once again how the interplay between matrix data and operator-algebraic structure governs the behaviour of the subfactors.


\newsection{Preliminaries}\label{Sec2}

We begin with briefly discussing two special classes of subfactors, known as the {\em Hadamard/spin model} and the {\em vertex model} subfactors. Their importance has been emphasized by Jones \cite{Jo2}.  We recall them briefly here, and interested readers can see \cite{JS} for details, or Section $2$ in \cite{BGG2} for a concise overview. 

Throughout the article, we shall use the shorthand notation $M_n$ to denote the set of all $n\times n$ complex matrices, and $\Delta_n \subset M_n$ denotes the unital subalgebra of the diagonal matrices in $M_n$. We always use the standard embedding of $M_n$ into $M_k \otimes M_n=M_{kn}$ by the map
\[
x \mapsto I_k \ot x \quad \text{for all }x \in M_n.
\]
\begin{dfn}
A complex Hadamard matrix is a unitary matrix in $M_n$ such that each entry is of modulus $1/\sqrt{n}$.
\end{dfn}

Let $F_n$ denote the discrete Fourier matrix $({\omega^{ij}}/{\sqrt{n}})_{i,j=0,1,...,n-1}$ of order $n$, where $\omega=e^{2\pi\sqrt{-1}/n}$ is a primitive $n$-th root of unity. Two complex Hadamard matrices $H_1, H_2 \in M_n$ are said to be {\em Hadamard equivalent}, denoted by $H_1\simeq H_2$, if there exist diagonal unitary matrices $D_1, D_2 \in \dt$ and permutation matrices $P_1, P_2 \in M_n$ such that  
\begin{equation}\label{he}
H_2=D_1 P_1 H_1P_2 D_2.
\end{equation}
This implements an equivalence relation on the set of all complex Hadamard matrices of order $n$. The {\em Hadamard equivalence class} of a complex Hadamard matrix $u$ is denoted by $[U]$. It is known that for $n = 2,3,5$, all the complex Hadamard matrices are Hadamard equivalent to the Fourier matrix $F_n$, whereas for $n=4$, there exists a one-parameter continuous family of complex Hadamard matrices up to equivalence. Note that for $n \geq 6$, a complete classification of complex Hadamard matrices up to Hadamard equivalence remains open in the literature.  

Complex Hadamard matrices are a rich source of hyperfinite subfactors, or more specifically, Hadamard/spin model subfactors. Any complex Hadamard matrix $U \in M_n$ induces the following
\[
\begin{matrix}
\Delta_n &\subset & M_n \cr \cup &\ & \cup\cr  \mathbb{C} &\subset &  \mathrm{Ad}_U(\Delta_n)
\end{matrix}
\]
non-degenerate commuting square (\cite{Po,JS}). Iterating Jones' basic construction on this non-degenerate commuting square horizontally, a subfactor $R_U \subset R$ is obtained, which is called the Hadamard subfactor. This subfactor is always irreducible, and the Jones index equals $n$. For more details, see \cite{JS}. Note that if $U$ lies in the Hadamard equivalence class of the Fourier matrix $F_n$, the subfactor $R_U\subset R$ is a crossed product by an outer action of $\bbz_n$ on $R$.

Now, let $U$ and $V$ be two distinct complex Hadamard matrices of the same order, and consider the corresponding Hadamard subfactors $R_U \subset R$ and $R_V \subset R$. Although $U \neq V$, the associated subfactors may coincide, that is, $R_U = R_V$, and we obtain a single Hadamard subfactor. To this end, we have a characterization result obtained in \cite{BG1}. Given $U\neq V$, we write $U \sim V$, if there exists a permutation matrix $P\in M_n$ and a diagonal unitary matrix $D\in\Delta_n$ such that $V=UPD$. This implements an equivalence relation '$\sim$' on the set of all complex Hadamard matrices of order $n$ that is finer than the Hadamard equivalence relation '$\simeq$' described in \Cref{he}. We have the following result.
\begin{thm}[Theorem $4.7$, \cite{BG1})]\label{main thm}
Let $U \neq V$ be complex Hadamard matrices of order $n$. Then, $R_U=R_V$ if and only if $U \sim V$. Therefore, if $U$ and $V$ are Hadamard inequivalent, that is, $ U \not\simeq V$, then the corresponding Hadamard subfactors $R_U\subset R$ and $R_V \subset R$ are always distinct.
\end{thm}

Another important class of hyperfinite subfactors is the so called {\em vertex model} subfactors.
\begin{dfn}[\cite{JS,BINA}]
A unitary matrix $U = (U^{\beta b}_{\alpha a}) \in M_n \otimes M_k$ is said to be a biunitary matrix if its block-transpose $W = (W^{\beta b}_{\alpha a})$, defined by $W^{\beta b}_{\alpha a}:= U^{\alpha b}_{\beta a}$, is also a unitary matrix in $M_n \otimes M_k$.
\end{dfn}

Let $W$ be a unitary matrix in $M_{nk}$. Consider the following quadruple of finite-dimensional algebras:
\[
\begin{matrix}
\bbc\otimes M_k  &\subset & M_n \otimes M_k \\
\cup & & \cup \\
\mathbb{C} &\subset & \mathrm{Ad}_W(M_n \otimes \mathbb{C})
\end{matrix}
\]
This forms a non-degenerate commuting square if and only if $W$ is a biunitary matrix \cite{JS}. Iterating Jones' basic construction horizontally, we obtain a subfactor $R_W \subset R$ with index $[R: R_W] = k^2$. This class of subfactors is called the vertex model subfactors. In contrast to the Hadamard subfactors, vertex model subfactors need not be irreducible, but they are always of integer index.



\newsection{Vertex model subfactors from complex Hadamard matrices}\label{Sec3}

It is a natural question whether there is any crucial relation between the Hadamard and the vertex model subfactors. In \cite{BGG2}, the authors investigated this question and obtained infinitely many vertex model subfactors (in fact, a nested chain of such subfactors) arising from a fixed {\em pair} of complex Hadamard matrices of the same order. All these vertex model subfactors have the Hadamard subfactor as an intermediate. We recall a few essential ingredients and results from that paper that we shall use here.

Let $U=(U_{ij})$ be a complex Hadamard matrix of order $n$ and consider the diagonal unitary matrix $D_U:=\sum_{i,j=1}^{n} \overline{U_{ij}} (E_{ii} \ot E_{jj})$ in $M_n\ot M_n=M_{n^2}$. 
For any complex Hadamard matrix $w \in M_n(\mathbb{C})$, consider the following finite-dimensional unital subalgebra
\[
A:= \text{Ad}_{U_1}(W \dt W^*) \subset \dt  \ot M_n \nonumber 
\]
where $U_1:=UD_U$. Denote the unique normalised tracial state on $M_n\ot M_n$ by $tr$. Let $E_1,E_2$ denote the unique trace preserving conditional expectations from $\dt  \ot M_n$ onto $A$ and $M_n$ respectively. Consider the following quadruple
\begin{IEEEeqnarray*}{lCl}\label{bi-unitary matrix}
\Gamma:=\quad\begin{matrix}
M_n &\subset^{\,E_2} & \dt \otimes M_n\\
\,\cup^{\,tr} & & \, \, \cup^{\,E_1}\\
\bbc &\subset^{\,tr} & A
\end{matrix}
\end{IEEEeqnarray*}
of finite-dimensional unital algebras. Then, $\Gamma$ is a nondegenerate commuting square (see Lemma 3.6 of \cite{BGG2} in this regard), that is, $E_1E_2=E_2E_1$ and $ A \cap M_n=\bbc$. Iterating Jones' basic construction on $\Gamma$, we obtain the following grid 
\begin{eqnarray}\label{construction of first vertex}
\begin{matrix}
M_n &\subset & \dt \ot M_n & \subset^{e_1 \ot I_n} & M_n \ot M_n  & \subset^{e_2 \ot I_n}  &\cdots & \subset R\\
\cup & & \cup & & \cup &  \, \, \, \, \, &&\cup\\
\bbc &\subset &  A &  \subset  &\langle A, e_1 \ot I_n \rangle &  \subset  & \cdots &  \subset R^{u,w}
\end{matrix}
\end{eqnarray}
of algebras, where $e_1=\frac{1}{n}\sum_{i,j}E_{ij} \in M_n$, $e_2=\sum_{i} E_{ii} \ot E_{ii} \in \dt \ot M_n$ etc. are the Jones' projections. The unital subalgebra
\begin{eqnarray}\label{first vertex}
R^{U,W}:=\overline{\bigvee\{A, \,e_1 \ot I^{(k)}_n, e_2 \ot I^{(k)}_n : k \in \bbn\}}^{sot} \subset R
\end{eqnarray}
is a subfactor of the hyperfinite type ${II}_1$ factor $R$. 
We recall the following fundamental result from \cite{BGG2}. 
In that work, the notation $R^{U,W}_0$ was used instead of $R^{U,W}$.

\begin{thm}[\cite{BGG2}, Theorem 4.1] \label{verttextospin}  
The inclusion $R^{U,W}\subset R$ is a vertex model subfactor of index $n^2$. Moreover, the Hadamard subfactor $R_U \subset R$ is an intermediate subfactor, i.e.,
\[
R^{U,W} \subset R_U \subset R.
\]
\end{thm}
Now, let $W$ be any complex Hadamard matrix, and $U$ be a complex Hadamard matrix in the Hadamard equivalence class $[W]$. By \Cref{he}, we write $u=DPW\widetilde{P}\widetilde{D}$ for some diagonal unitary matrices $D,\widetilde{D}\in\Delta_n$ and permutation matrices $P,\widetilde{P}\in M_n$. Further, applying \Cref{main thm}, we can discard $\widetilde{D},\widetilde{P}$ and assume that $U=DPW$. We recall the following result.

\begin{ppsn}[Theorem $5.2$, \cite{BGG2}]\label{proposition}
Let $U, V \in[W]$ be Hadamard equivalent complex Hadamard matrices and $R_U, R_V \subset R$ be the corresponding Hadamard subfactors. Assume that $U=D_1P_1W$ and $V=D_2P_2W$. For any arbitrary complex Hadamard matrix $W$ of the same order as $U$ and $V$, there is an automorphism $\theta\in\mathrm{Aut}(R)$ such that $\theta(R_U)=R_V$ and the subfactors $R^{U,W} \subset R\mbox{ and }R^{V,W}\subset R$ are isomorphic. Moreover,
\begin{enumerate}[$(i)$]
\item if $P_1 \ne P_2$, then $\theta\in\mathrm{Out}(R)$;
\item if $P_1 = P_2$, then $\theta= \mathrm{Ad}_{D_2 D_1^*} \in \mathrm{Inn}(R)$.
\end{enumerate}
\end{ppsn}

Condition $(ii)$ above precisely says that the subfactors $R_U ,R_V \subset R$ are conjugate to each other when $P_1=P_2$, and vice-versa. Note that Theorem $5.2$ in \cite{BGG2} is a much stronger result; however, we only require the above in this paper. In fact, it is possible to construct the automorphism $\theta$ explicitly. To see this, we adopt the notation $M_n^{(k)}$ for the $k$-fold tensor product of $M_n$ with itself, and write $A^{(k)}$ for $A^{\otimes k}$ when $A \in M_n$. Let $D_1, D_2 \in \dt$ be diagonal unitary matrices and $P_1, P_2 \in M_n$ be permutation matrices as defined in \Cref{proposition}, and set $P:=P_2P_1^*$. Then, the sequence $\{\mathrm{Ad}_{P^{(k)}} : k \in \mathbb{N}\}\subseteq\mathrm{Inn}(R)$ converges pointwise to an automorphism $\widetilde{\theta} \in \operatorname{Aut}(R)$ in the strong operator topolgy (see Proposition $5.4$, \cite{BGG2}). The automorphism $\theta$ is then given by the composition $\mathrm{Ad}_{D_2 D_1^*} \circ \widetilde{\theta}$. Note that if $P_1=P_2$, then $\widetilde{\theta}=id_{R}$, and consequently, $\theta = \mathrm{Ad}_{D_2 D_1^*}$.


\newsection{Main results}\label{Sec4}

In this section, we prove Theorem A and B stated in the Introduction.

\subsection{Factoriality of the intersection}

Let us begin by fixing some notations to be used frequently throughout the rest of the article.
\begin{notation}\label{matrices}\rm
Let $\vec{n} = (n_1, n_2, \ldots, n_k) \in \mathbb{N}^k$ be a fixed vector with each $n_i \geq 2$.
\begin{enumerate}[$(i)$]
\item For a complex Hadamard matrix  $U=(U_{ij})$ order $n$, define $D_U :=\sum_{i,j=1}^{n} \overline{U_{ij}} (E_{ii} \ot E_{jj})$, a diagonal unitary matrix in $M_n \ot M_n$. Also, consider the unitary matrix $U_1:=U D_U \in \Delta_n \ot M_n$.
\item $F_{n_i}$ denotes the Fourier matrix of order $n_i$, and consider the $k-$fold tensor product
\begin{equation}\label{W}
W:=F_{n_1} \otimes F_{n_2} \otimes \cdots \otimes F_{n_k}.
\end{equation}
This matrix will be used throughout the paper repeatedly.
\item For each $n \in \mathbb{N}$ and $k \in \{0, 1, \ldots, n-1\}$, define the diagonal matrices:
    \begin{align*}
    \mathscr{D}_{n,0} := I_n, \quad
    \mathscr{D}_{n,1} := \mathrm{diag}(1, \omega, \omega^2, \ldots, \omega^{n-1}), \quad
    \mathscr{D}_{n,k} := \mathscr{D}_{n,1}^k = \mathrm{diag}(1, \omega^k, \omega^{2k}, \ldots, \omega^{(n-1)k}),
    \end{align*}
where $\omega = e^{-2\pi\sqrt{-1}/ n}$ is a primitive $n$-th root of unity.
\item Define the cyclic permutation matrices:
    \begin{align*}
    \sigma_{n,0} := I_n, \quad \quad
    \sigma_{n,1} := E_{12} + E_{23} + \cdots + E_{(n-1)n} + E_{n1}, \quad \quad
    \sigma_{n,k} := \sigma_{n,1}^k.
    \end{align*}
\item For any vector $\vec{r} = (r_1, r_2, \ldots, r_k) \in \mathbb{N}^k$, define
\[
    \mathscr{D}_{\vec{r}} := \mathscr{D}_{n_1, r_1} \otimes \mathscr{D}_{n_2, r_2} \otimes \cdots \otimes \mathscr{D}_{n_k, r_k}, \quad
    \sigma_{\vec{r}} := \sigma_{n_1, r_1} \otimes \sigma_{n_2, r_2} \otimes \cdots \otimes \sigma_{n_k, r_k}.
\]
\end{enumerate}
\end{notation}

\noindent The following matrix identities are easy to verify\,:
\begin{align}\label{first}
\mathscr{D}_{n,1} F_n &= F_n \sigma_{n,n-1}, &
\sigma_{n,1} F_n &= F_n \mathscr{D}_{n,1},
\end{align}
Therefore, 
\begin{align}\label{second}
\mathrm{Ad}_{F_n^*}\mathscr{D}_{n,k} &= \sigma_{n,n-k}, &
\mathrm{Ad}_{F_n}\mathscr{D}_{n,k} &= \sigma_{n,k}.
\end{align}

We start by establishing a sequence of results needed to prove Theorem A.
\begin{lmma}\label{lemma1}
Let $\vec{n} = (n_1, n_2, \ldots, n_k)$ be the fixed vector described as in \Cref{matrices}. For any vector $\vec{r} = (r_1, r_2, \ldots, r_k) \in \mathbb{N}^k$ with $1 \leq r_i \leq n_i$ for all $i$, let $\mathscr{D}_{\vec{r}}$ denote the diagonal unitary matrix and $\sigma_{\vec{r}}$ the corresponding permutation matrix as defined in \Cref{matrices}. Let $W:= F_{n_1} \otimes F_{n_2} \otimes \cdots \otimes F_{n_k}$ be the $k$-fold tensor product of the Fourier matrices $F_{n_i}$ of order $n_i$. The following identities hold\,:
\begin{enumerate}[$(i)$]
    \item $\mathrm{Ad}_{W}(\mathscr{D}_{\vec{r}}) = \sigma_{\vec{r}}$,
    \item $\mathrm{Ad}_{W^*}(\mathscr{D}_{\vec{r}}) = \sigma_{\vec{n} - \vec{r}}$,
\end{enumerate}
where $\vec{n} - \vec{r} := (n_1 - r_1, n_2 - r_2, \ldots, n_k - r_k)$.    
\end{lmma}
\begin{prf} From the identities in \Cref{second}, we have $\mathrm{Ad}_{F_{n_i}}(\mathscr{D}_{n_i, r_i}) = \sigma_{n_i, r_i}$ for each $i = 1, \ldots, k$. Since $W = F_{n_1} \otimes \cdots \otimes F_{n_k}$ and $\mathscr{D}_{\vec{r}} = \bigotimes_{i=1}^k \mathscr{D}_{n_i, r_i}$, it follows that
\[
\mathrm{Ad}_{W}(\mathscr{D}_{\vec{r}}) = \bigotimes_{i=1}^k \mathrm{Ad}_{F_{n_i}}(\mathscr{D}_{n_i, r_i}) = \bigotimes_{i=1}^k \sigma_{n_i, r_i} = \sigma_{\vec{r}}\,.
\]
Similarly, since $W^* = F_{n_1}^* \otimes \cdots \otimes F_{n_k}^*$ and $\mathrm{Ad}_{F_{n_i}^*}(\mathscr{D}_{n_i, r_i}) = \sigma_{n_i, n_i - r_i}$, we obtain
\[
\mathrm{Ad}_{W^*}(\mathscr{D}_{\vec{r}}) = \bigotimes_{i=1}^k \sigma_{n_i, n_i - r_i} = \sigma_{\vec{n} - \vec{r}}\,.
\]
This completes the proof.\qed
\end{prf}

The following result is in the spirit of \cite{Bur}.
\begin{lmma}\label{action of group 1}
Let $W$ be the complex Hadamard matrix $F_{n_1} \otimes F_{n_2} \otimes \cdots \otimes F_{n_k}$ and $R_U \subset R$ be the Hadamard subfactor corresponding to the complex Hadamard matrix $U \in [W]$. Without loss of generality, in the spirit of \Cref{main thm}, assume that $U=DPF$ for some diagonal unitary matrix $D\in\Delta_N$ and permutation matrix $P\in M_N$. Then, 
\[
R = R_{U} \rtimes_{\theta} G.
\]  
where $\theta$ is an outer action of the abelian group $G=\mathbb{Z}_{n_1} \times \mathbb{Z}_{n_2} \times ... \times \mathbb{Z}_{n_k}$ on $R_U$ defined by  
\[
\theta_{\Vec{r}} = \text{Ad}_{ DP\mathscr{D}_{\Vec{r}}P^*D^*}
\]  
for all $\Vec{r}=(r_1, r_2,...,r_k) \in G$.
\end{lmma}
\begin{prf} From the construction of Hadamard subfactor $R_U \subset R$ (recall from \Cref{Sec2}), we have the following 
\begin{eqnarray}\label{spin dfn}
R_U&=& \overline{\bigvee \{\mathrm{Ad}_{U}(  \Delta_N ), \,e_2 \ot I^{(k-1)}_N, \,e_1 \ot I^{(k)}_N :  k \in \bbn \} }^{sot}\subset R
\end{eqnarray}
Now, we claim that for $\Vec{r} \in G$, we have $\theta_{\Vec{r}}(R_U)=R_U$. First, observe the following
\begin{IEEEeqnarray*}{lCl}
\theta^\prime_{\Vec{r}}(U\Delta_N U^*)&=&\text{Ad}_{DP\mathscr{D}_{\Vec{r}} P^*D^*}(U \Delta_N U^*) \nonumber \\
 &=&\text{Ad}_{\{ DP \mathscr{D}_{\Vec{r}}\}}( W \Delta_N W^*) \nonumber \hspace{2mm} \\
     &=&\text{Ad}_{\{DP  \mathscr{D}_{\Vec{r}} \}}((W\Delta_N W^*)) \nonumber \hspace{3mm}  \\
    &=&\text{Ad}_{U_1}(\mathscr{D}_{\Vec{r}} W \Delta_N W^* \mathscr{D}^*_{\Vec{r}}) \nonumber \quad \quad (\text{by \Cref{lemma1}})\\
    &=&\text{Ad}_{DP}(W \sigma_{\Vec{r}} \Delta_N \sigma^*_{\Vec{r}} W^*) \nonumber \\ 
   &=&\text{Ad}_{DP}(W \Delta_N  W^*) \nonumber\\
   &=&U\Delta_N U^*.
\end{IEEEeqnarray*}
Since for $i \in \bbn$, the Jones' projections $e_1 \ot I^{(i)}_N$ and $e_2 \ot I^{(i-1)}_N$ commutes with $U \mathscr{D}_{\Vec{r}} U^*$ for all $\Vec{r} \in G$, it follows that
\begin{eqnarray}\label{fix projection}
\theta_{\Vec{r}}(e_1 \ot I^{(i)}_N)=e_1 \ot I^{(i)}_N \quad  \text{and} \quad \theta_{\Vec{r}}(e_2 \ot I^{(i-1)}_n)=e_2 \ot I^{(i-1)}_N.
\end{eqnarray}
By the continuity of $\theta_{\Vec{r}}$ and \Cref{spin dfn} to \Cref{fix projection}, we conclude that
\[
\theta_{\Vec{r}}(R_U) = R_U.
\]  
Hence, $\theta_{\Vec{r}}$ is well-defined for all $\Vec{r} \in G$, which completes the proof of the claim. 

By the construction of Hadamard subfactor $R_U \subset R$, we have the following non-degenerate commuting square. 
\[
\begin{matrix}
R_{U} &\subset & R\\
\cup & & \cup\\
\bbc &\subset & \Delta_N
\end{matrix}
\]
and in particular, $R_U\cap\Delta_N=\bbc$ and $R_U\Delta_N=\Delta_NR_U=R$. We show that $\theta_{\Vec{r}} \in \operatorname{Out}(R_U)$ for all $\Vec{r} \neq 0 \in G$. Let $y \in R_U$ satisfies $\theta_{\Vec{r}}(x)y=yx$ for all $x \in R_U$, and then $(\text{Ad}_{DP}(\mathscr{D}_{\Vec{r}}))^{*}y \in R\cap R_{U}^{'}$. Since $R_U\subset R$ is irreducible, it follows that $y =\alpha \text{Ad}_{DP}(\mathscr{D}_{\Vec{r}})$ for some $\alpha\in \bbc$. Suppose that $\alpha \neq 0$, which implies $\mathrm{Ad}_{DP}(\mathscr{D}_{\Vec{r}}) \in R_U  \cap \Delta_N$. Since $R_U \cap \Delta_N =\bbc$, it easily follows that $\text{Ad}_{DP}(\mathscr{D}_{\Vec{r}}) \in \bbc$. This contradicts the fact that $\mathscr{D}_{\Vec{r}} \not \in \bbc$, $\Vec{r} \neq 0 \in G$. Hence $\alpha=0$, which then implies that $y=0$. This concludes that $\theta_{\Vec{r}} \in \operatorname{Out}(R_U)$ for all $\Vec{r} \neq 0 \in G$. Since $R_U \Delta_N=R$ and $\{\text{Ad}_{DP}(\mathscr{D}_{\Vec{r}}) : \Vec{r} \in G\}$ form a basis of $\Delta_N$, we conclude that $R=R_{U}\rtimes^{\theta} G$, which completes the proof. \qed
\end{prf}

Suppose that $M$ be an type $II_1$ factor and $\theta: G \to \mathrm{Aut}(M)$ be an outer action of a finite group $G$ on $M$, i.e., $\theta_g\in\mathrm{Out}(M)$ for each $g\neq 1$. In this situation, the associated crossed product $M \rtimes_\theta G$ is again a type $II_1$ factor. It is known that $M\subset M\rtimes G$ is a finite index subfactor with $[M\rtimes G:M]=|G|$. We also have the fixed-point (unital) subalgebra
\[
M^G := \big\{ x \in M : \theta(g)(x) = x \text{ for all } g\in G \big\} \subseteq M.
\]
Recall, $M \subset M \rtimes_\theta G$ is irreducible, and so is $M^G \subset M$. Moreover, $(M^G)' \cap (M\rtimes_{\theta} G)=\mathbb{C}G$, where $\mathbb{C}G\subseteq M\rtimes G$ denotes the group algebra of $G$.
\begin{ppsn}\label{group}
Let $\theta : G \to \operatorname{Aut}(M)$ be an outer action of a finite group $G$ on the type $II_1$ factor $M$. Then, for any unitary element $u \in \mathbb{C}G$, the intersection $\mathrm{Ad}_u(M) \cap M $ is always a type $II_1$ factor. Moreover, the Jones index $[M:\mathrm{Ad}_u(M) \cap M]$ is finite and equal to the order of some subgroup $H\leq G$.
\end{ppsn}
\begin{prf}
Let $u \in \mathbb{C}G\subset M\rtimes G$ be a unitary element. Observe that, $u$ normalizes the fixed-point algebra $M^G$, that is $\mathrm{Ad}_u(M^G) = M^G$. In particular, this implies $M^G = \mathrm{Ad}_u(M^G) \subset \mathrm{Ad}_u(M)$. Therefore, we get the following
\[
M^G \subset \, M \cap \mathrm{Ad}_u(M) \, \subset M.
\]
Since the fixed-point subalgebra $M^G \subset M$ is irreducible, it immediately follows that the intersection $M \cap \mathrm{Ad}_u(M)$ is a subfactor of $M$, and $M \cap \mathrm{Ad}_u(M)=M^H$ for some subgroup $ H$ of $G$. Hence, we have the index $[M: \mathrm{Ad}_u(M) \cap M]=|H|$.\qed
\end{prf}

We are now ready to prove our first main result.
\medskip

\noindent\textbf{\em Proof of Theorem A:\,} Since the Hadamard subfactors are conjugate to each other, due to \Cref{proposition}, we can assume that $U = D P W$ and $V = \widetilde{D} P W$, where $W=F_{n_1} \otimes F_{n_2} \otimes F_{n_3} \otimes\ldots\otimes F_{n_k}$, for some diagonal unitary matrices $D, \widetilde{D}$ and permutation matrix $P$ of order $N=:n_1n_2\cdots n_k$.  Again from \Cref{proposition}, it follows that the corresponding Hadamard subfactors $R_U \subset R$ and $R_V \subset R$ are isomorphic under the inner automorphism $\theta = \mathrm{Ad}_{\widetilde{D}D^*} \in \mathrm{Inn}(R)$. That is,  
\begin{equation}\label{conjugate spin}
\theta(R_U) = \mathrm{Ad}_{\widetilde{D}D^*}(R_U) = R_V.
\end{equation}
By \Cref{action of group 1}, we have $R = R_U \rtimes_{\alpha} G$ for an outer action $\alpha$ of a finite abelian group $G$, whose group algebra is $\mathbb{C} G = \Delta_N \subset R$. Setting $M = R_U$ and $u = \widetilde{D}D^*$, where $u$ is a unitary element of $\mathbb{C} G = \Delta_N$, in \Cref{group}, we finally conclude that
\[
R_U \cap \mathrm{Ad}_{\widetilde{D}D^*}(R_U) = R_U \cap R_V
\]
(by \Cref{conjugate spin}) is a factor.\qed
\smallskip

\textbf{Question:} Does Theorem A holds if we remove the hypothesis that the Hadamard subfactors are conjugate?


\subsection{Index and relative commutant}

Now that we know that $R_U\cap R_V\subset R$ is always a subfactor when $R_U,R_V\subset R$ are conjugate to each other and $U,V \in [W]$, it is natural to investigate the index $[R:R_U \cap R_V]$ and the relative commutant $(R_U \cap R_V)^\prime\cap R$. First we establish a few technical results.
 
\begin{lmma}\label{lemma2}
Let $ k \in\bbn$ and $\vec{n} = (n_1, n_2, \ldots, n_k)\in\mathbb{N}^k$ be a fixed vector. Consider the complex Hadamard matrix $V_k=F_{n_1} \otimes F_{n_2} \otimes F_{n_k}$ of order  $N_k=n_1\ldots n_k$, where $F_{n_j}$'s are the Fourier matrices of order $n_j$. Then, we have $(V_k)_1=P_kW$, where $P_k$ is a block-diagonal permutation matrix of order $N_k^2$ and $(V_k)_1=V_kD_{V_k}$ (see \Cref{matrices}).
\end{lmma}
\begin{prf}
 For any $n \in \bbn$, by definition, we have $(F_n)_1=F_{n} D_{F_n}$ (see \Cref{matrices})  and by using \Cref{second}, we get $F_{n} D_{F_n} = P F_{n}$ where $P := \sum_{i=1}^{n_1} E_{ii} \otimes \sigma_{n,n - i + 1}$ is a permutation matrix (see \Cref{matrices}). Therefore, it follows that $(F_n)_1=PF_n$; this concludes the base step.

Assume that the result holds for $(t-1)-$th step. That is, $(V_{t-1})_1=V_{t-1} D_{V_{t-1}} = P_{t-1} V_{t-1}$ for some block-diagonal permutation matrix $P_{t-1}$ of order $N_{t-1} = n_1 \cdots n_{t-1}$. From the definition of $V_t$ and $D_{V_t}$ (see \Cref{matrices}) we have
\[
V_t = V_{t-1} \otimes F_{n_t}, \quad \text{and} \quad D_{V_t} = \sum_{i=1}^{n_t} E_{ii} \otimes D_{V_{t-1}} \otimes \mathscr{D}_{n_t,i}^*.
\]
Therefore, we get the following 
\begin{align*}
(V_t)_1 &=V_t D_{V_t}  \\
&= (V_{t-1} \otimes F_{n_t}) \left( \sum_{i=1}^{n_t} E_{ii} \otimes D_{V_{t-1}} \otimes \mathscr{D}_{n_t,i}^* \right) \\
&= \sum_{i=1}^{n_t} E_{ii} \otimes V_{t-1} D_{V_{t-1}} \otimes F_{n_t} \mathscr{D}_{n_t,i}^* \quad \text{(by induction hypothesis and \Cref{second})} \\
&= \sum_{i=1}^{n_t} E_{ii} \otimes P_{t-1} V_{t-1} \otimes \sigma_{n_t,i}^*  \\
&= \left( \sum_{i=1}^{n_t} E_{ii} \otimes P_{t-1} \otimes \sigma_{n_t,i}^* \right)(V_{t-1} \otimes F_{n_t}) \\
&= P_t V_t,
\end{align*}
where $P_t := \sum_{i=1}^{n_t} E_{ii} \otimes P_{t-1} \otimes \sigma_{n_t,i}^*$ is a permutation matrix of order $N_t$. This completes the proof.\qed
\end{prf}

\begin{ppsn}\label{lemma3}
Consider the Hadamard subfactor $R_U \subset R$ corresponding to the complex Hadamard matrix $U$ of order $N$ and $W=F_{n_1} \otimes F_{n_2} \otimes \cdots \otimes F_{n_k}$ where $n_1\ldots n_k=N$. Obtain the vertex model subfactor $R^{U,W}\subset R$ described in \Cref{Sec3}. Then, $(R^{U,W})^\prime\cap R_U=\bbc$ and $(R^{U,W})^\prime\cap R=\mathbb{C}^N$.
\end{ppsn}
\begin{prf}
Recall that the vertex model subfactor $R^{W,W}\subset R$ is obtained by iterating Jones' basic construction on the following non-degenerate commuting square
\[
\Gamma:=\quad\begin{matrix}
M_N & \subset & \Delta_N \otimes M_N \\
\cup & & \cup \\
\mathbb{C} & \subset & \mathrm{Ad}_{W_1}(W \Delta_N W^*)
\end{matrix}
\]
as described in \Cref{construction of first vertex}. By the Ocneanu compactness theorem (see Theorem 5.7.1, \cite{JS}), it follows that the relative commutant is given by
\begin{equation}\label{Ocneanu}
(R^{W,W})'\cap R=\left( \mathrm{Ad}_{W_1}(W \Delta_N W^*) \right)' \cap M_N,
\end{equation}
where $W_1 = W D_W \in \Delta_N \otimes M_N$ (see \Cref{matrices}).

By \Cref{lemma2}, we know that $W_1:=WD_W=PW$ for some block-diagonal permutation matrix $P \in \Delta_N \otimes M_N$. Since $\mathrm{Ad}_{F_n}(\Delta_n)=\mathrm{Ad}_{F^*_n}(\Delta_n)$ for all $n \in \bbn$, it follows that $\mathrm{Ad}_W (W \Delta_N W^*) = \Delta_N$. Therefore, we get the following
\begin{align*}
\mathrm{Ad}_{W_1}(W \Delta_N W^*) 
&= \mathrm{Ad}_{P W}(W \Delta_N W^*)  \\
&= \mathrm{Ad}_{P} \circ \mathrm{Ad}_W (W \Delta_N W^*) \\
&= \mathrm{Ad}_P (\Delta_N).
\end{align*}
Substituting this into \Cref{Ocneanu}, and using the fact $P$ is a block-diagonal permutation matrix of order $N^2$, we obtain
\begin{equation}\label{commutant of vertex}
(R^{W,W})'\cap R = \left( \mathrm{Ad}_P (\Delta_N) \right)' \cap M_N = \mathbb{C}^N.
\end{equation}
Further, by the construction of Hadamard subfactors, we know that
\begin{eqnarray*}
\begin{matrix}
\Delta_N  &\subset & R\\
\cup & & \cup\\
\bbc &\subset &  R_W
\end{matrix}
\end{eqnarray*}
is a non-degenerate commuting square. Therefore, it follows that $R_W \cap \Delta_N=\bbc$. Hence, from \Cref{commutant of vertex}, we get that $(R^{W,W})^\prime\cap R_W= R_W \cap \bbc^N=\bbc$.

We have the Hadamard subfactor $R_U\subset R$ associated to the complex Hadamard matrix $U\in [W]$. Without loss of generality, we can assume that $U=DPW$ due to \Cref{main thm}, for some diagonal unitary matrix $D$ and permutation matrix $P$, each of order $N=n_1n_2\cdots n_k$. By \Cref{proposition}, we know that the following chain of subfactors 
\[
R^{W,W}\subset R_W\subset R \quad \mbox{and} \quad R^{U,W}\subset R_U \subset R
\]
are isomorphic. Therefore, by \Cref{commutant of vertex} and the fact that $(R^{W,W})^\prime\cap R_W=\bbc$, we conclude that
\begin{equation} \label{relative commutant}
(R^{U,W})^\prime\cap R = \mathbb{C}^N
\end{equation}
and $(R^{U,W})^\prime\cap R_U=\bbc$.\qed
\end{prf}

\begin{ppsn}\label{action of group}
Let $W$ be the complex Hadamard matrix $F_{n_1} \otimes F_{n_2} \otimes \cdots \otimes F_{n_k}$ and $R_U \subset R$ be the Hadamard subfactor corresponding to the complex Hadamard matrix $U \in [W]$. Then, 
\[
R_U = R^{U,W}\rtimes_{\theta'} G
\]  
where $\theta'$ is an outer action of the abelian group $G:=\mathbb{Z}_{n_1} \times \mathbb{Z}_{n_2} \times ... \times \mathbb{Z}_{n_k}$ on the type $II_1$ factor $R^{U,W}$ defined by  
\[
\theta'_{\Vec{r}} = \text{Ad}_{(U \mathscr{D}_{\Vec{r}} U^*)}
\]  
for all $\Vec{r}=(r_1, r_2,...,r_k) \in G$.
\end{ppsn}
\begin{prf}
By \Cref{first vertex}, we have the following vertex model subfactor 
\begin{eqnarray}\label{vertex of fourier}
R^{U,W} &=& \overline{\bigvee \{\mathrm{Ad}_{U_1}( W \Delta_N W^*), \,e_2 \ot I^{(k)}_N, \,e_1 \ot I^{(k)}_N :  k \in \bbn \} }^{sot}\subset R
\end{eqnarray}
We claim that for $\Vec{r} \in G$, $\theta^{'}_{\Vec{r}}(R^{U,W})=R^{U,W}$. 
First, observe the following
\begin{IEEEeqnarray*}{lCl}
\theta^{'}_{\Vec{r}}(\text{Ad}_{U_1}(W \Delta_N W^*)) &=& \text{Ad}_{U_1 \mathscr{D}_{\Vec{r}}  U_1^*}(\text{Ad}_{U_1}( W \Delta_N W^*))\hfill{\qquad(\mbox{since }U_1=(I_n \ot U) D_U)}\\
&=&\text{Ad}_{U_1 ( \mathscr{D}_{\Vec{r}} ) U^*_1}(\text{Ad}_{U_1}( F_n\dt F^{*}_n)\\
&=&\text{Ad}_{U_1}(\mathscr{D}_{\Vec{r}} W \Delta_N W^* \mathscr{D}^*_{\Vec{r}})\hfill{(\text{by \Cref{second}})}\\
&=&\text{Ad}_{U_1}(W \sigma_{\Vec{r}} \Delta_N \sigma^*_{\Vec{r}} W^*)\\ 
&=&\text{Ad}_{U_1}(W \Delta_N  W^*).
\end{IEEEeqnarray*}
Since for $i \in \bbn$, the Jones' projections $e_1 \ot I^{(i)}_N$ and $e_2 \ot I^{(i)}_N$ commutes with $U \mathscr{D}_{\Vec{r}} U^*$ for all $\Vec{r} \in G$, it follows that
\begin{eqnarray}\label{equation4}
\theta_{\Vec{r}}^{'}(e_1 \ot I^{(i)}_N)=e_1 \ot I^{(i)}_N \quad  \text{and} \quad \theta_{\Vec{r}}^{'}(e_2 \ot I^{(i)}_n)=e_2 \ot I^{(i)}_N.
\end{eqnarray}
Hence, by the continuity of $\theta^{'}_{\Vec{r}}$ and using \Cref{vertex of fourier} to \Cref{equation4}, we conclude that
\[
\theta^{'}_{\Vec{r}}(R^{U,W})=R^{U,W}.
\]  
Therefore, $\theta^{'}_{\Vec{r}}$ is well-defined for all $\Vec{r} \in G$, which concludes the proof of the claim. 

By the construction of the Hadamard subfactor $R_U \subset R$, we have the following non-degenerate commuting squares
\[
\begin{matrix}
M_N &\subset & R\\
\cup & & \cup\\
U \Delta_N U^* &\subset &   R_{U},
\end{matrix}\qquad\mbox{ and }\qquad\begin{matrix}
 M_N  &\subset & R\\
\cup & & \cup\\
\bbc &\subset &R^{U,W}
\end{matrix}
\]
(see \Cref{construction of first vertex} of \cite{BGG2}), and therefore it follows that the following quadruple 
\[
\begin{matrix}
U\Delta_N U^* &\subset & R_U\\
\cup & & \cup\\
\bbc &\subset &  R^{U,W}
\end{matrix}
\]
is also a non-degenerate commuting square. Hence, we have 
\begin{eqnarray}\label{large intersection}
R^{U,W}\cap \mathrm{Ad}_U( \Delta_N) =\bbc   \quad \text{and} \quad R^{U,W}(\mathrm{Ad}_{U}(\Delta_N))=R_U.
\end{eqnarray}
We now aim to prove that $\theta'$ is an outer action of the group $G$. Specifically, we claim that $\theta'_{\Vec{r}} \in \text{Out}(R^{U,W})$ for all $\Vec{r}\neq \Vec{0}\in G$, where $\Vec{0}$ is the identity element of $G$. Consider any $y \in R^{U,W}$ satisfying $\theta'_{\Vec{r}}(x)y = yx$ for all $x \in R^{U,W}$. This implies that $(\text{Ad}_{U}(\mathscr{D}_{\Vec{r}}))^* y$ lies in $(R^{U,W})'\cap R_U$. Therefore, by \Cref{relative commutant} it follows that $y$ must be of the form $y = \alpha \, \text{Ad}_{U}(\mathscr{D}_{\Vec{r}})$ for some $\alpha \in \bbc$. Suppose that $\alpha \neq 0$. Then, it follows that $\text{Ad}_{U}(\mathscr{D}_{\Vec{r}})$ lies in $\mathrm{Ad}_U(\Delta_N) \cap R^{U,W}$. Using \Cref{large intersection}, we have $\mathrm{Ad}_U(\Delta_N) \cap R^{U,W}=\bbc$. Therefore, $\text{Ad}_U(\mathscr{D}_{\Vec{r}}) \in \bbc$, which implies that $\mathscr{D}_{\Vec{r}} \in \bbc$. However, this leads to a contradiction since $\mathscr{D}_{\Vec{r}} \not \in \bbc$ for $\Vec{r} \neq \Vec{0} \in  G$. Hence, we must have $ y = 0$. This proves that $\theta^{'}_{\Vec{r}}$ is an outer automorphism on the vertex model subfactor $R^{U,W}\subset R$. Furthermore, since $\mathscr{D}_{\Vec{0}} = I_N$, it follows that $\theta^{'}_{\Vec{0}}$ is the identity automorphism. Therefore, we conclude that $\theta'$ is an outer action on the vertex model subfactor $R^{U, W}\subset R$. 

Finally, observe that the following  set  
\[
\{ \text{Ad}_U(\mathscr{D}_{\Vec{r}}) :  \Vec{r} \in G \}  
\]  
forms a basis for $U \Delta_N U^*$ (see \Cref{matrices}). Using \Cref{large intersection}, we conclude that $R_U = R^{U, W}\rtimes_{\theta'} G $, which completes the proof.\qed
\end{prf}

Recall that $\Delta_N\subset M_N$ denotes the unital subalgebra of diagonal matrices. Given complex Hadamard matrices $U, V \in [W]$, consider the following abelian unital subalgebra of $M_N$\,:
\begin{eqnarray}\label{finite}
\mathscr{A}_U^V:=\mathrm{Ad}_U(\Delta_N)\cap\mathrm{Ad}_V(\Delta_N).
\end{eqnarray}
Clearly, $\dim(\mathscr{A}_U^V)\leq N=\dim(\Delta_N)$. We are now ready to prove our second main result.
\medskip

\noindent\textbf{\em Proof of Theorem B:\,} By \Cref{action of group}, we have $R_U = R^{U,W}\rtimes_{\theta'} G$. Therefore, $R_U$ can be expressed as
\begin{equation*}
R_U = \left\{ \sum_{\Vec{r}} x_{\Vec{r}} \left( U \mathscr{D}_{\Vec{r}} U^* \right) : \, x_{\Vec{r}} \in R^{U,W}, \, \Vec{r} \in G \right\}.
\end{equation*}
Since $R^{U,W} \subset R_U \cap R_V\subset R_U$, it follows that
\[
R_U \cap R_V = R^{U,W} \rtimes_{\theta'} H
\]
where $H$ is a subgroup of $G$. Hence $[R_U \cap R_V : R^{U,W}]=|H|$, and consequently, we have
\begin{equation}\label{intersection set}
R_U \cap R_V = \left\{ \sum_{\Vec{r}} x_{\Vec{r}} \left( U \mathscr{D}_{\Vec{r}} U^* \right) : x_{\Vec{r}} \in R^{U,W}, \,  \Vec{r} \in H \right\}.
\end{equation}
By the multiplicativity of the Jones index, we obtain that $[R:R_U\cap R_V]=N^{2}/|H|$.
 
Now, we show that $|H| = \dim(\mathscr{A}_U^V)$. By construction of the vertex model subfactor $R^{U,W} \subset R$, we have the following non-degenerate commuting square
\[
\begin{matrix}
M_N & \subset & R \\
\cup & & \cup \\
\mathbb{C} & \subset & R^{U,W}
\end{matrix}
\]
Hence, for any $x \in R^{U,W}$ and $\mathrm{Ad}_{U}({D}) \in U \Delta_N U^*$, we have
\begin{equation*}
E^R_{M_N}(x(\mathrm{Ad}_U(D))) = E_{M_N}(x) \, \mathrm{Ad}_U (D)  = \tau(x) \, \mathrm{Ad}_{U}(D),
\end{equation*}
where $\tau$ denotes the unique trace on the hyperfinite type $II_1$ factor $R$. Applying this to \Cref{intersection set}, we obtain
\begin{equation}\label{expectation}
E_{M_N}(R_U \cap R_V) = \left\{ \sum_{\Vec{r}} \alpha_{\Vec{r}} \left( U \mathscr{D}_{\Vec{r}} U^* \right) : \alpha_{\Vec{r}} \in \mathbb{C},  \Vec{r} \in H\right\}.
\end{equation}
Since both the quadruples $(U\Delta_N U^* \subset M_N, \, R_U \subset R)$ and $(V\Delta_NV^* \subset M_N, \, R_V \subset R)$ are non-degenerate commuting squares, it follows that the following quadruple
\[
\begin{matrix}
R_U \cap R_V & \subset & R \\
\cup & & \cup \\
\mathscr{A}_U^V  & \subset & M_N
\end{matrix}
\]
is also a commuting square. Hence, we conclude that $E^R_{M_N}(R_U \cap R_V)=\mathscr{A}_U^V$. Substituting this into \Cref{expectation}, we obtain the following
\begin{equation}\label{1st intersection}
\mathscr{A}_U^V=\left\{\sum_{\Vec{r}} \alpha_{\Vec{r}} \left( U \mathscr{D}_{\Vec{r}} U^* \right) : \alpha_{\Vec{r}} \in \mathbb{C}, \,\Vec{r} \in H \right\}.
\end{equation}
Therefore, the set $\{ U \mathscr{D}_{\Vec{r}} U^* : \Vec{r} \in H \}$ forms a basis for $\mathscr{A}_U^V$, and hence $\dim(\mathscr{A}_U^V) = |H|$. This completes the proof of the claim.

To complete the proof, it is enough to show that $(R_U \cap R_V)' \cap R=(\mathscr{A}_U^V)'\cap\bbc^N$. By \Cref{large intersection} and \Cref{1st intersection}, it follows that 
\begin{equation}\label{AB}
R_U \cap R_V = R^{U,W}\mathscr{A}_U^V
\end{equation}
and hence, we have the following
\[
(R_U \cap R_V)' \cap R = \left( (R^{U,W})' \cap R \right) \cap ((\mathscr{A}_U^V)' \cap R).
\]
By \Cref{relative commutant}, we have $(R^{U,W})' \cap R = \mathbb{C}^N$. Substituting this in the above equation, we finally obtain
\begin{eqnarray*}
(R_U \cap R_V)' \cap R &=& \mathbb{C}^N \cap ((\mathscr{A}_U^V)' \cap R) \nonumber \\
                       &=& \mathbb{C}^N \cap(\mathscr{A}_U^V)',
\end{eqnarray*}
which concludes the proof.\qed

\begin{ppsn}\label{attain}
For every subgroup $H \leq G$, there exists a pair of Hadamard subfactors $R_U,R_V\subset R$ corresponding to the complex Hadamard matrices $U,V\in [F_{n_1}\otimes\ldots\otimes F_{n_k}]$, such that
\[
[R : R_U \cap R_V]=N^2/|H|
\]
where $N=n_1n_2 \cdots n_k$.
\end{ppsn}
\begin{prf}
We begin by establishing the following\,:
\smallskip

\noindent\textbf{Claim:} For any divisor $k$ of $n \in \mathbb{N}$, there exists a pair of complex Hadamard matrices $u, v \in [F_n]$ such that $\mathrm{dim}(\mathrm{Ad}_u\Delta_n)\cap\mathrm{Ad}_v(\Delta_n))=k$, where $v=Du$ for some diagonal unitary matrix $D\in\Delta_n$.

Let $l \in \mathbb{N}$ such that $n = k l$. Define the diagonal unitary matrix $D \in \Delta_n$ by
\[
D = \sum_{r=1}^{k} \sum_{s=1}^{l} \alpha_r E_{s,s},
\]
where $\alpha_1, \dots, \alpha_k$ are distinct complex numbers of modulus one, chosen in such a way that $D$ is not a scalar multiple of diagonal unitary matrices $\mathscr{D}_{n,m}$ for $m \in \{1, 2, \ldots, n\}$ (see \Cref{matrices}). Note that otherwise, we would have $v = F_n \sigma_{n,m}$ (see \Cref{second}), and hence, $\mathrm{dim}(u \Delta_n u^* \cap v \Delta_n v^*) = n$. We have
\[
\mathrm{Ad}_v(\mathscr{D}_{n,i})=\beta_i \mathrm{Ad}_u(\mathscr{D}_{n,i})
\]
if and only if $i \in \{l, 2l, \ldots, kl\}$. Therefore,
\[
u \Delta_n u^* \cap v \Delta_n v^*=\big\{\sum_{i=1}^{k} r_i \mathrm{Ad}_{u}(\mathscr{D}_{n,il}) \; : \; r_i \in \mathbb{C}\big\}.
\]
Hence, $\dim(u \Delta_n u^* \cap v \Delta_n v^*) = k$, which completes the proof of the claim.

Now, let $H$ be a subgroup of the finite abelian group $G = \mathbb{Z}_{n_1} \times \cdots \times \mathbb{Z}_{n_k}$. Then $H = \mathbb{Z}_{m_1} \times \cdots \times \mathbb{Z}_{m_k}$ where each $m_i$ divides $n_i$. From the claim above, for each $i = 1, 2, \ldots, k$, there exists a diagonal unitary matrix $D_i$ of order $n_i$ such that
\begin{equation} \label{last eqn}
\dim\big(\mathrm{Ad}_{F_{n_i}}(\Delta_{n_i}) \cap \mathrm{Ad}_{D_i F_{n_i}}(\Delta_{n_i})\big)=m_i.
\end{equation}
Let $U = W = F_{n_1} \otimes \cdots \otimes F_{n_k}$ and define $D = D_1 \otimes \cdots \otimes D_k$ so that $V = D U$ is also a complex Hadamard matrix of order $N = n_1 n_2 \cdots n_k$. Using \Cref{last eqn}, we conclude the following\,:
\[
\dim(U \Delta_N U^* \cap V \Delta_N V^*) = m_1 m_2 \cdots m_k = |H|.
\]
Finally, from Theorem C, it follows that 
\[
[R : R_U \cap R_V]=N^2/|H|,
\]
which completes the proof.\qed
\end{prf}

\subsection{Intersection as a vertex model subfactor}

Recall the construction of the {\em vertex model} subfactors from \Cref{Sec2}, and the abelian unital subalgebra $\mathscr{A}_U^V$ of $M_N$ from \Cref{finite}.

\begin{thm}\label{intersection is vertex}
The intersection $R_U \cap R_V \subset R$ is a vertex model subfactor if and only if $\mathscr{A}_U^V=\mathbb{C}$.
\end{thm}
\begin{prf}
By \Cref{AB}, we have $R_U \cap R_V = R^{U,W}\mathscr{A}_U^V$. If $\mathscr{A}_U^V=\mathbb{C}$, then clearly $R_U \cap R_V\subset R$ is a vertex model subfactor. For the converse direction, we show that whenever $\mathrm{dim}(\mathscr{A}_U^V)>1,\,R_U \cap R_V\subset R$ is not a vertex model subfactor. Recall from \Cref{Sec2} that a vertex model subfactor $R_Z \subset R$, where $Z\in M_N \otimes M_N$ is a bi-unitary matrix, is constructed by iterating the basic construction on the following
\[
    \begin{matrix}
    \bbc\otimes M_N & \subset & M_N \otimes M_N \\
    \cup & & \cup \\
    \mathbb{C} & \subset & \mathrm{Ad}_Z(M_N \otimes \mathbb{C})
    \end{matrix}
\]
non-degenerate commuting square. Therefore,
\begin{equation}\label{holds}
R_Z\cap M_N =\mathbb{C}. 
\end{equation}
However, from the definition of $\mathscr{A}_U^V$, we have $\bbc\subsetneq\mathscr{A}_U^V\subseteq M_N$ and $\mathscr{A}_U^V\subseteq R_U \cap R_V$. Therefore, $(R_U \cap R_V) \cap M_N \neq \mathbb{C}$. Hence, the intersection fails to satisfy \Cref{holds}. Thus, $R_U \cap R_V \subset R$ is not a vertex model subfactor.\qed
\end{prf}

The above result generalizes our earlier works in \cite{BG1,BGG}, where the orders of the complex Hadamard matrices considered were small.


\newsection{Application to the Connes-St{\o}rmer relative entropy}\label{Sec5}

In this final section, we analyse the relative position of two conjugate Hadamard subfactors through an invariant. The {\em Connes--St{\o}rmer relative entropy} provides a quantitative measure of the mutual position of two subalgebras of a finite von Neumann algebra. For inclusions of the form $R_U, R_V \subset R$, this entropy captures how far the two subfactors are from being simultaneously contained in a common intermediate subfactor.

Let $\mathcal{M}$ be a finite von Neumann algebra equipped with a fixed faithful normal trace. For finite-dimensional subalgebras $\mathcal{P},\mathcal{Q}\subset\mathcal{M}$, Connes \& St{\o}rmer introduced the notion of relative entropy $H(\mathcal{P|Q})$ between $\mathcal{P}\mbox{ and }\mathcal{Q}$. In \cite{PP}, Pimsner \& Popa observed that the definition of $H(\mathcal{P|Q})$ works as it is for arbitrary von Neumann subalgebras $\mathcal{P},\mathcal{Q}\subset\mathcal{M}$. For more on entropy, see \cite{NS}.

In our situation, we have a pair of subfactors $R_U, R_V \subset R$, and in the present subsection, we determine legitimate lower and upper bounds for the Connes--St{\o}rmer relative entropy $H(R_U|R_V)$. Note that the quadruple $(R_U\cap R_V\subset R_U, R_V\subset R)$ often fails to become a commuting square. Therefore, our situation is a bit more complicated since results from \cite{WW} is not applicable here.

\begin{dfn}[\cite{CS}]
Let $\mathcal{P, Q } \subseteq  \mathcal{M}$ be von Neumann subalgebras of the finite von Neumann algebra $(\mathcal{M},\tau)$. Let $S$ denotes the set of finite families $\gamma=(x_1, x_2,...,x_n)$ of positive elements in $\mathcal{M}$ with $\sum_{i=1}^nx_i=1$. Then,
\[
H_\gamma(\mathcal{P|Q}):= \sum_{j=1}^n\left(\tau\circ\eta\,E_{\mathcal{Q}}^{\mathcal{M}}(x_j)-\tau\circ\eta\,E_{\mathcal{P}}^{\mathcal{M}}(x_j)\right)
\]
where $\eta :[0,\infty)\rightarrow\mathbb{R}$ is the continuous function $t\mapsto -t\log t$. The Connes-St{\o}rmer relative entropy between $\mathcal{P}\mbox{ and }\mathcal{Q}$ is defined as $H(\mathcal{P|Q}):=\sup_{\gamma \in S}\,H_\gamma(\mathcal{P|Q})$.
\end{dfn}

It is known that finding the exact value of the relative entropy is often very hard. In \cite{BG2}, the authors have calculated the exact value of the relative entropy for a very special instance. In the literature, a modified version of the Connes-St{\o}rmer relative entropy exists due to Choda \cite{choda2} (see also \cite{choda}), which is less complicated to compute. This is defined as $h(\mathcal{P}|\mathcal{Q}):=\sup_{\gamma \in S}\,h_\gamma(\mathcal{P|Q})$, where
\[
h_\gamma(\mathcal{P|Q}):= \sum_{j=1}^n\left(\tau\circ\eta\,E_{\mathcal{Q}}^{\mathcal{M}}(x_j)-\tau\circ\eta\,(x_j)\right).
\]

We recall the following result from \cite{BG1}.

\begin{ppsn}[Proposition $2.6$, \cite{BG1}]\label{popaadaptation2}
\begin{enumerate}[$(i)$]
\item Let $\{M_n\}, \{A_n\}$ and $\{B_n\}$ be increasing sequences of von Neumann subalgebras of a finite von Neumann algebra $M$ such that $\{A_n\}, \{B_n\}\subset M$ and $M=\big(\bigcup_{n=1}^{\infty}M_n\big)^{\dprime}.$ If $A=\big(\bigcup_{n=1}^{\infty} A_n\big)^{\dprime}$ and $B=\big(\bigcup_{n=1}^{\infty} B_n\big)^{\dprime}$, then $H(B|A)\leq\displaystyle \liminf_{n \to \infty}H(B_n|A_n).$
\item If in addition, $E_{A_{n+1}}E_{M_n}=E_{A_n}$ and $E_{B_{n+1}}E_{M_n}=E_{B_n}$ for $n\in \mathbb{N}$, then $H(B|A)=\displaystyle\lim_{n \to \infty}H(B_n|A_n)$ increasingly.
\end{enumerate}
\end{ppsn}

Note that a similar statement also holds for `$h$' in place of `$H$'. Our main result in this section is the following.

\begin{thm}\label{ent of cs}
Let $U,V\in [F_{n_1}\otimes\ldots\otimes F_{n_k}]$ be distinct complex Hadamard matrices and consider the corresponding Hadamard subfactors $R_U\subset R$ and $R_V\subset R$. If the Hadamard subfactors are conjugate, then we have the following\,:
\begin{enumerate}[$(i)$]
\item The modified Connes-St{\o}rmer entropy is given by 
\[
h(R_U|R_V)=\frac{1}{N}\sum_{i,j= 1}^{N} \eta |(U^{*}V)_{ij}|^{2},
\]
where $(U^{*}V)_{ij}$ denotes the $ij-$th entry of the matrix $U^{*}V$ .
\item One has the following\,:
\begin{equation*}
h(R_U|R_V)\leq H(R_U|R_V)\leq\log\big(N/\mathrm{dim}(\mathscr{A}_U^V)\big)
\end{equation*}
where $\mathscr{A}_U^V=\mathrm{Ad}_U(\Delta_N)\cap\mathrm{Ad}_V(\Delta_N)$ as defined in \Cref{finite}.
\end{enumerate}
\end{thm}
\begin{prf}
Since $R_U \subset R$ and $R_V \subset R$ are conjugate subfactors by  applying \Cref{proposition}, we assume that $R_U = \mathrm{Ad}_D(R_{PW})$ and $R_V = \mathrm{Ad}_{\widetilde{D}}(R_{PW})$, where $R_{PW} \subset R$ denotes the Hadamard subfactor associated with the complex Hadamard matrix $PW$  and $D, \widetilde{D}$ are diagonal unitary matrices. Therefore, we have the following
\begin{eqnarray}\label{entropy equation}
h(R_U \,|\, R_V) &=& h\big(\mathrm{Ad}_{D}(R_{PW}) \,\big|\, \mathrm{Ad}_{\widetilde{D}}(R_{PW})\big) \nonumber\\
&=& h \big(R_{PW} \,\big|\, \mathrm{Ad}_{D^{*}\widetilde{D}}(R_{PW})\big).
\end{eqnarray}
From \Cref{action of group 1}, it is known that $R = R_{PW} \rtimes^{\theta} G$, where the outer action $\theta$ is given by  $\theta_{\vec{r}} = \mathrm{Ad}_{\left(P\, \mathscr{D}_{\vec{r}}\, P^{*}\right)}$ for all $ \vec{r} \in G$. Since $\{\mathscr{D}_{\vec{r}} : \vec{r} \in G\}$ forms a basis for the diagonal subalgebra $\Delta_{N} \subset M_N$, where $N = n_{1} n_{2} \cdots n_{k}$, a diagonal unitary matrix $D^{*}\widetilde{D}$ can be written as  
\[
D^{*}\widetilde{D}=\sum_{i=0}^{n-1} \tau\big(D^{*}\widetilde{D} \, P^{*} \mathscr{D}^{*}_{\vec{r}} P\big) \, P^{*} \mathscr{D}_{\vec{r}} P.
\]
and therefore, applying Theorem $3.14$ of \cite{choda2}, from \Cref{entropy equation} we obtain the following
\begin{eqnarray}\label{lower bound}
h(R_U \,|\, R_V) 
&=& h\big(R_{PW} \,\big|\, \mathrm{Ad}_{D^{*}\widetilde{D}}(R_{PW})\big) \nonumber\\
&\leq& \sum_{\vec{r} \in G} \eta\big(\tau\big(D^{*}\widetilde{D} \, P^{*} \mathscr{D}_{\vec{r}} P\big) \; \tau\big(D^{*}\widetilde{D} \, P^{*} \mathscr{D}^{*}_{\vec{r}} P\big) \big) \nonumber\\
&=& \sum_{\vec{r} \in G} \eta\,\big|\tau\big(D^{*}\widetilde{D} \, P^{*} \mathscr{D}_{\vec{r}} P\big)\big|^{2}.
\end{eqnarray}
By Corollary $3$ of \cite{choda}, we have the following  
\begin{eqnarray}\label{upper bound}
h(U\Delta_NU^*\,\big|\,V\Delta_NV^*) 
&=& h\big(\Delta_N\,\big|\,U^*V\Delta_NV^*U\big) \nonumber\\
&=& \frac{1}{N} \sum_{i,j=1}^{N} \eta\big(\big|(U^{*}V)_{ij}\big|^{2} \big).
\end{eqnarray}
Now to compute the $(i,j)$-th entry of $U^{*}V$, observe that  
\begin{eqnarray}\label{matrix entries}
U^{*}V  &=& \mathrm{Ad}_{W^{*}}\big(P D^{*}\widetilde{D} P\big) \nonumber\\
&=& \mathrm{Ad}_{W^{*}}\Big(\sum_{\vec{r} \in G} \tau\big( P D^{*}\widetilde{D} P^{*} \mathscr{D}^{*}_{\vec{r}} \big) \, \mathscr{D}_{\vec{r}} \Big)\hfill{\qquad(\text{apply \Cref{lemma1}})}\nonumber\\
&=& \sum_{\vec{r} \in G} \tau\big( P D^{*}\widetilde{D} P^{*} \mathscr{D}^{*}_{\vec{r}} \big) \, \sigma_{\vec{n} - \vec{r}}
\end{eqnarray}
Since each $\sigma_{\vec{n} - \vec{r}}$ is a permutation matrix of order $N$ (see \Cref{matrices}), it follows from \Cref{matrix entries} that for every $\vec{r} \in G$, the scalar $\tau(P D^{*} \widetilde{D} P^{*} \mathscr{D}^{*}_{\vec{r}})$ occurs exactly $N$-times among the entries of the unitary matrix $U^{*}V$. Moreover, from the definition of the permutation matrices $\{ \sigma_{\vec{r}}: \vec{r} \in G \}$, we conclude that for each $\vec{r}\in G$ that there exists unique $j \in \{1,2,..,N\}$ such that $(U^{*}V)_{1j}= \tau\big( P D^{*}\widetilde{D} P^{*} \mathscr{D}^{*}_{\vec{r}} \big)$. Hence, it follows that  
\begin{eqnarray}\label{upper bound of mass}
\sum_{j=1}^{N}\eta |(U^*V)_{1j}|^2 
&=&\sum_{\vec{r} \in G} \eta |\tau(D^*\widetilde{D})|^{2} \nonumber \\
&=&\frac{1}{N}\sum_{i,j= 1}^{N} \eta |(U^{*}V)_{ij}|^{2}
\end{eqnarray}
By \Cref{popaadaptation2}, we know that  $h(U\Delta_NU^*\big|V\Delta_NV^*)\leq h(R_U|R_V)$. Combining this with \Cref{lower bound,upper bound of mass}, we deduce the following
\[
h(R_U|R_V)=\frac{1}{N}\sum_{i,j= 1}^{N} \eta |(U^{*}V)_{ij}|^{2},
\]
which establishes the first part of the theorem.
\smallskip

Next, by \Cref{lemma3} we have $(R^{U,W})^\prime\cap R_U=\bbc$. Since $R^{U,W}\subset R_U \cap R_V$, it follows that $(R_U \cap R_V)^\prime\cap R_U=\bbc$. Applying Corollary $4.6$ of \cite{PP} and using Theorem C, we obtain the following\,:
\[
H(R_U \,|\, R_U \cap R_V)=\log([R_U : R_U \cap R_V])=\log\big(N/\dim(\mathscr{A}_U^V)\big).
\]
Moreover, by definition we have $h(R_U|R_V)\leq H(R_U|R_V)$, and from \cite{PP} it follows that $H(R_U|R_V)\leq H(R_U|R_U\cap R_V)$. Combining these inequalities, we finally get the following  
\[
h(R_U|R_V)\leq H(R_U|R_V)\leq\log\left(N/\dim(\mathscr{A}_U^V)\right),
\]
which finishes the proof.\qed
\end{prf}

{\em Concluding remarks:\,} Note that in general, $H(R_U|R_V)\leq H(R|R_V)=\log[R:R_V]=\log N$. Therefore, the upper bound $\log\big(N/\mathrm{dim}(\mathscr{A}_U^V)\big)$ is indeed an improvement for the relative entropy. The interplay between index, relative commutant, and entropy gives a complete picture of the intersection of conjugate Hadamard subfactors. The index $[R:R_U \cap R_V]$ measures the algebraic size of the intersection, the relative commutant records its internal symmetry, and the entropy $H(R_U|R_V)$ captures its analytic separation.  
Together, these invariants reveal a rich geometric structure underlying the space of Hadamard subfactors associated with tensor products of Fourier matrices. Furthermore, the explicit formulae obtained here suggest further directions\,: extending the analysis to non-conjugate pairs, more general Hadamard matrices which are not necessarily from the equivalence classes of DFT-type matrices, and a pair of biunitary matrices.

\bigskip

\bigskip

\noindent {\em Department of Mathematics and Statistics},\\
{\em Indian Institute of Technology Kanpur},\\
{\em Uttar Pradesh $208016$, India}
\medskip

\noindent {\bf Keshab Chandra Bakshi:} keshab@iitk.ac.in, bakshi209@gmail.com\\
{\bf Satyajit Guin:} sguin@iitk.ac.in\\
{\bf Guruprasad:} guruprasad8909@gmail.com


\bigskip

\begin{thebibliography}{10}

\bibitem{BG1}
Bakshi, K. C.; Guin, S.:
\newblock Relative position between a pair of spin model subfactors.
\newblock {\em J. Aust. Math. Soc.} 119 (2025), no. 1, 1--38.

\bibitem{BG2}
Bakshi, K. C.; Guin, S.:
\newblock A family of subfactors arising from a pair of complex Hadamard matrices.
\newblock {\em Internat. J. Math.} 36 (2025), no. 2, Paper No. 2450071.

\bibitem{BGG}
Bakshi, K. C.; Guin, S.; Gururprasad:
\newblock Conjugate pairs of Hadamard subfactors and vertex models.
\newblock {\em Proc. Indian Acad. Sci. Math. Sci.} 135 (2025), no. 2, Paper No. 23.

\bibitem{BGG2}
Bakshi, K. C.; Guin, S.; Gururprasad:
\newblock A noncommutative construction of families of biunitary matrices and application to subfactors.
\newblock {\em arXiv:2506.14351}.

\bibitem{BINA}
Bhattacharyya, B.: 
\newblock Group actions on graphs related to Krishnan-Sunder subfactors.
\newblock {\em Trans. Amer. Math. Soc.} 355 (2003), no. 2, 433--463.

\bibitem{Bur}
Burstein, R.:
\newblock Group-type subfactors and Hadamard matrices.
\newblock {\em Trans. Amer. Math. Soc.} 367 (2015), no. 10, 6783--6807.

\bibitem{choda}
Choda, M:
\newblock Relative entropy for maximal abelian subalgebras of matrices and the entropy of unistochastic matrices.
\newblock{\em Internat. J. Math.} 19 (2008), no. 7, 767--776.

\bibitem{choda2}
Choda, M:
\newblock Conjugate pairs of subfactors and entropy for automorphisms.
\newblock{\em Internat. J. Math.} 22 (2011), no. 4, 577--592.

\bibitem{CS}
Connes, A.; Størmer, E.:
\newblock Entropy for automorphisms of $II_1$ von Neumann algebras.
\newblock {\em Acta Math.} 134 (1975), no. 3-4, 289--306.

\bibitem{Connes}
Connes, A.:
\newblock Periodic automorphism of the hyperfinite factor of type $II_1$.
\newblock {\em Acta Sci. Math.} 39 (1977), no. 1-2, 39--66.

\bibitem{Jo}
Jones, V. F. R.:
\newblock Index for subfactors.
\newblock {\em Invent. Math.} {\bf 72} (1983), no. 1, 1--25.

\bibitem{JS}
Jones, V. F. R.; Sunder, V. S.:
\newblock Introduction to subfactors.
\newblock London Mathematical Society Lecture Note Series, 234. {\em Cambridge University Press, Cambridge}, 1997.

\bibitem{JX}
Jones, V. F. R.; Xu, F.:
\newblock Intersections of finite families of finite index subfactors.
\newblock {\em Internat. J. Math.} 15 (2004), no. 7, 717--733.

\bibitem{Jo1}
Jones, V. F. R.:
\newblock Two subfactors and the algebraic decomposition of bimodules over $II_1$ factors.
\newblock {\em Acta Math. Vietnam.} 33 (2008), no. 3, 209--218.

\bibitem{Jo2}
Jones, V. F. R.:
\newblock Planar algebra, I.
\newblock{\em New Zealand J. Math.} 52, 1--107, 2021.

\bibitem{NS}
Neshveyev, S; Størmer, E.:
\newblock Dynamical entropy in operator algebras.
\newblock Springer-Verlag, Berlin, 2006.

\bibitem{PP}
Pimsner, M.; Popa, S.:
\newblock Entropy and index for subfactors.
\newblock {\em Ann. Sci. École Norm. Sup.} (4) 19 (1986), no. 1, 57--106.

\bibitem{Po}
Popa, S.:
\newblock Classification of amenable subfactors of type II.
\newblock {\em Acta Math.} 172 (1994), no. 2, 163--255.

\bibitem{SW}
Sano, T.; Watatani, Y.:
\newblock Angles between two subfactors.
\newblock {\em J. Operator Theory} 32 (1994), no. 2, 209--241.

\bibitem{WW}
Watatani, Y.; Wierzbicki, J.:
\newblock Commuting squares and relative entropy for two subfactors.
\newblock {\em J. Funct. Anal.} 133 (1995), no. 2, 329--341.

\end{thebibliography}
\end{document}